\documentclass[a4paper,10pt]{article}
\pdfoutput=1
\usepackage[T1]{fontenc}
\usepackage[utf8]{inputenc}
\usepackage{lmodern}
\usepackage{hyperref}
\usepackage{graphicx}
\usepackage[english]{babel}
\usepackage{xcolor}
\usepackage{url}

\usepackage{bbm}

\usepackage{float}
\usepackage{amsthm}
\usepackage{amsmath}
\usepackage{amssymb}
\usepackage{graphicx}
\usepackage{setspace}
\usepackage{esint}
\usepackage{hyperref}
\usepackage{comment}
\usepackage{ulem}
\usepackage{enumitem}
\usepackage{rotating}
\usepackage{enumitem}
\usepackage{subcaption}
\usepackage{pdfpages}

\newtheorem{theorem}{Theorem}

\newcommand{\EE}{\mathbb{E}}
\newcommand{\PP}{\mathbb{P}}
\newcommand{\RR}{\mathbb{R}}

\newcommand{\EU}{\Phi}

\theoremstyle{definition}

\renewcommand{\emph}{\textit}

\RequirePackage[includefoot]{geometry}
\title{\Large\bf The interdependence between hospital choice and waiting time --- with a case study in urban China}

\author{\normalsize
Joris van de Klundert, Business School Universidad Adolfo Ib\'a\~nez 
\and  \normalsize
Roberto Cominetti, Faculty of Engineering, Universidad Adolfo Ib\'a\~nez 
\and \normalsize
Yun Liu, School of Health Policy and Management, Erasmus University Rotterdam
\and\normalsize
Qingxia Kong, Rotterdam School of Management, Erasmus University Rotterdam \\}
\date{}

\begin{document}

\maketitle

%
%
%
%
%

{\normalsize

\noindent {\bf Abstract:}
 Hospital choice models often employ random utility theory and include waiting time as a choice determinant. When applied to evaluate health system improvement interventions, these models disregard that hospital choice in turn is a determinant of waiting time. We present a novel, general model capturing the endogeneous relationship between waiting time and hospital choice, including the choice to opt out, and characterize the unique equilibrium solution of the resulting convex problem. 
 
We apply the general model in a case study on the urban Chinese health system, specifying that patient choice follows a multinomial logit (MNL) model and waiting times are determined by M/M/1 queues. The results reveal that analyses which solely rely on MNL models overestimate the effectiveness of present policy interventions and that this effectiveness is limited. We explore alternative, more effective, improvement interventions.\\


\noindent {\bf Funding:} The third author was supported by the Chinese Scholarship Council [grant number 201507720036] 
}

\vfill
\newpage 

\section{Introduction}

Long waiting presents a significant barrier to healthcare access in various forms and in a variety of health services and settings across the globe \cite{mcintyre2020waiting}. It often disadvantages patients of lower socioeconomic status and is negatively associated with patient outcomes such as mortality \cite{mcintyre2020waiting, martinez2019prolonged, jones2022association}. 

For patients, waiting times are an important determinant of their choice of hospital and hospital level and may even cause them to avoid care and opt out \cite{victoor2012determinants, fischer2015understanding, liu2018factors}. The significance of waiting time as a determinant of hospital choice is evidenced by revealed preferences and by stated preferences obtained in a number of discrete choice experiments (DCEs) conducted across the globe \cite{brown2015hospital, berhane2015patients, zhu2019exploring, liu2020impact} (however  not for all patient populations and contexts \cite{smith2018patient,liu2018patients}). 

 Long waiting can be viewed to signal a need for health system interventions that address underlying performance issues. In urban China, for instance, long waiting has led to great discontent and caused undesired avoidance of care \cite{liu2020impact}. In this setting, long waiting times especially arise at overutilized tertiary hospitals, as patients freely access hospitals at the level of their choice and tend to leave primary care underutilized \cite{liu2020impact, arsenault2020hospital}. Both policy measures and modeling studies reported in the academic literature have sought to improve the persistent utilization problems by strengthening primary care \cite{li2017primary, liu2020impact, healthreform2020}. The model studies regard waiting time as a hospital attribute, or as a hospital level attribute, i.e., as an independent explanatory variable determining patient choice. 

As observed by Lee and Cohen, waiting times in hospitals may not be static, independent, hospital (level) attributes but typically vary with the choices made \cite{lee1985equilibrium, marianov2008facility, zhang2019medical}. Increased choice probabilities can cause congestion and result in longer waiting. Thus, the variable waiting time, which is presumed to be independent in the aforementioned DCE models, is endogeneously related to the dependent variable hospital choice. This form of endogeneity, i.e. simultaneity involving a dependent and a presumed independent is essentially different from endogeneity between independents, or between independents and unobserved factors, for which econometric methods are readily available (see \cite{train2009discrete} and the references therein).

Disregard of endogeneity in DCE models can imply inconsistency and incorrectness in the model estimation and cause the model to be ill suited for policy analysis \cite{train2009discrete, antonakis2010making}. In the case of hospital choice in urban China, it may explain why evidence based interventions to improve medical equipment and skills of physicians in primary care have not resolved the utilization problems \cite{meng2019can, ta2020trends, li2017primary, li2020quality, liu2018factors}. For example, interventions to strengthen primary care can increase the probability of choosing primary care and reduce the probability to choose tertiary care. The revised choice probabilities translate to a tertiary care waiting time decrease that in turn increases its choice probability and thus may partially undo the intervention effects. A similar effect may occur when interventions to reduce the opt out probability translate into longer waiting times of hospitals, partially undoing the relative increase in attractiveness. 

In this study, we advance the theoretical analysis of the endogenous relationship between waiting time and choice probability and present a case study analysing policy interventions in urban China (\cite{marianov2008facility, zhang2019medical,kucukyazici2020incorporating}). The presented general model incorporates the interdependence between waiting times and hospital choice, including the choice to opt out. It employs random utility theory to specify functions capturing how choice probabilities depend on utilities and how these utilities subsequently depend on waiting times. In addition, the model relies on a fairly general function to express how waiting times depend on choice probabilities. We show that the model yields unique equilibrium choice probabilities and waiting times, which can be characterized and found by convex optimization.

In the case study, we specifically apply a model in which a multinomial logit (MNL) model specifies choice probability as a function of waiting time (among other factors) and a queuing model specifies waiting time as a function of choice probability \cite{marianov2008facility,zhang2019medical}. This special case model resolves the disregard of endogeneity in previously published preference and policy studies that consider waiting time to be exogenous. It applies the presented model to provide a consistent assessment of the quantitative effects of health system interventions (such as interventions to strengthen primary care) 
that is more accurate than assessments obtained by previous MNL based models\cite{liu2020impact, zhu2019exploring}. The analysis reports on intervention effects obtained with the  
special case model and on significance of differences in policy intervention effects obtained through models presented in previous studies.

In the discussion, we extensively consider future research on hospital choice behaviour and models. 
	
\section{Models and Equilibrium Solutions}

The model deals with a finite set of healthcare facility types $i\in I=\{1,\ldots,n\}$, each one characterized by a set of attributes such as skills of medical staff, equipment, price of service, and capacity. The types can refer to hospital levels, such as primary, secondary, and tertiary hospitals. 

On the other hand, we consider a finite set $K=\{1,\ldots,m\}$ of patient types where each $k\in K$ represents a homogeneous patient population. The homogeneity refers to patient characteristics such as symptoms, insurance type, age, and occupation, that determine the utility $U^k_{i}$ attached to facility type $i$ by patients of type $k$. If the locations of facilities and patients are relevant choice determinants --- as might be the case when considering travel time --- each patient and/or facility may constitute a separate type by itself.

We model utilities in terms of random variables. Namely, for each patient type $k \in K$ and facility type $i \in I$, the utility is described by a
random variable $U^k_{i}=u^k_{i}+\varepsilon^k_{i}$ where $u^k_{i}\in\RR$ is
the expected utility level, and $\varepsilon^k_{i}$ is a random deviation 
that captures the variability within the $k$-th patient population with $\EE(\varepsilon^k_{i})=0$. The model also includes the possibility that patients choose to opt out. To this purpose, we introduce a fictitious opt-out facility type $i=0$ with utility $U^k_{0}=u^k_{0}+\varepsilon^k_{0}$ for patients of type $k \in K$ \cite{kucukyazici2020incorporating}. 

Following random utility theory, we assume patients choose the facility that yields the maximal utility and hence the probabilities $\pi^k=(\pi^k_{i})_{i\in I}$ of choosing facility type $i$ are given by the discrete choice model
$$\pi^k_{i}=\PP(U^k_{i}\geq U^k_{j},\,\forall j\in I\cup\{0\}).$$ 
For simplicity we assume that the random vector $\varepsilon^k=(\varepsilon^k_{i})_{i\in I}$ has a continuous distribution. 

Denoting $u^k=(u^k_{0},u^k_{1},\ldots,u^k_{n})$ and considering the expected utility function 
$$\EU^{k}(u^k)=\EE(\max\{u^k_{0}+\varepsilon^k_{0},u^k_{1}+\varepsilon^k_{1},\ldots,u^k_{n}+\varepsilon^k_{n}\}),$$ we explicitly establish that:\\

\begin{theorem}
$\EU^{k}(\cdot)$ is smooth and convex, and its derivatives are precisely the choice probabilities, namely
\begin{equation}\label{gradprox}
\pi^k_i=\pi^k_i(u^k)=\frac{\partial\EU^{k}}{\partial u^k_{i}}(u^k)\qquad(\forall i\in I\cup\{0\}).\\
\end{equation}
\end{theorem}

The proof is presented in Appendix A. We refer to \cite{domencich1975urban, sheffi1978another, williams1977formation, lee1985equilibrium} for related (partial) results.\\

In the remainder, we adopt the common assumption that utility varies linearly with the waiting time and that an increase (decrease) in waiting time translates into a decrease (increase) in utility. More formally, we consider utilities of the form
\begin{equation}\label{utility}
u^k_i=\bar u^k_i-\alpha^k w_i
\end{equation}
where $\bar u^k_i$ is a reference utility, $w_i$ is the expected waiting time at facility $i\in I$, and the parameter $\alpha^k>0$ represents the sensitivity of patients of type $k \in K$ to waiting time. We reflect on this common assumption in the discussion section as a topic for future research.\\

Let us now turn to expressing how choice probabilities determine waiting time. To this purpose, we first define for each patient type $k\in K$ the arrival rate $I^k>0$ as the number of  
patients of type $k$ arriving per time unit. Then, for facility type $i \in I$ and choice probabilities $\pi^k_i, k \in K$, and assuming Poisson arrivals, the total arrival rate at facility type $i$ can be expressed as
\begin{equation}
\label{flow} x_i  =  \sum_{k\in K}I^k \pi^k_i \qquad (\forall i\in I).
\end{equation}
The waiting time is assumed to be zero when opting out, i.e.,  $w_0\equiv 0$, whereas for $i\geq 1$ we assume that it is
given by a strictly increasing continuous function $\theta_i:[0,\bar x_i)\to (0,\infty)$ of the arrival rates, namely 
\begin{equation}
\label{wait} w_i = \theta_i(x_i)\qquad\hspace{3.3ex} (\forall i\in I),\\
\end{equation}
with $\bar x_i>0$ a saturation level at which waiting times diverge, that is $$\lim_{x\to\bar
x_i}\theta_i(x_i)=\infty.$$ 
The saturation level $\bar x_i$ can serve to model the finite capacity of each of the facility types and the infinitely long queues that can form when patient volumes are at or above capacity. However, the model allows for $\bar x_i=\infty$ at some or all facility types (as in \cite{lee1985equilibrium}).

Given the exogenous rates  $I^k>0$ and the parameters $\bar u^k_i$ and $\alpha^k$ in the utilities \eqref{utility}, we look for $w_i$ 
and $x_i$, $i\in I$, simultaneously satisfying the system of equations \eqref{gradprox}-\eqref{wait}. The equations \eqref{flow} and \eqref{wait} extend the standard equations \eqref{gradprox}, \eqref{utility} in random utility based choice models, thus capturing the endogenous relationship between choice and waiting time.\\

\noindent{\sc Remark.}
From \eqref{flow} it follows that any solution satisfies $x_i\geq 0$ and therefore $w_i\geq w^0_i$ 
where $ w_i^0\triangleq \theta_i(0)$. We observe that the equalities $x_i=0$ and $w_i= w^0_i$ can 
not be excluded at equilibrium. 
This will occur for a facility type $i \in I$ whenever the choice 
probabilities are zero for all $k\in K$, that is $\pi_i^k=0$, which is the case if the random terms 
$\varepsilon_i^k$ have bounded support and the facilities $j\neq i$ provide much higher utilities so that $i$ is never optimal.\\

Let us consider a continuous and strictly increasing extension of $\theta_i(\cdot)$ to the negative reals with $\theta_i(x_i)\to-\infty$ for $x_i\to-\infty$. Since any equilibrium
satisfies $x_i\geq 0$, the specific form of the extension 
 of $\theta_i(\cdot)$ is irrelevant, and we may simply set $\theta_i(x_i)=\theta_i(0)+x_i$ for $x_i<0$. This extended function has a well defined inverse $\theta_i^{-1}:(-\infty,\infty)\to(-\infty,\bar x_i)$ which is also strictly increasing and continuous. Thus, we can eliminate $x_i$ from  \eqref{flow}-\eqref{wait} and express these equations in terms of the variables $w_i$ alone, namely
\begin{equation}\label{reduced}
\theta_i^{-1}(w_i)=\sum_{k\in K}I^k\pi^k_i \qquad (\forall i\in I).
\end{equation}

The equilibrium choice probabilities and waiting times can now be characterized as follows: \\

\begin{theorem} \label{thm1} The system of equations \eqref{gradprox}-\eqref{wait} has a unique solution. This equilibrium is the unique minimizer of the strictly convex smooth function $\Theta:\RR^d\to\RR$ given by
\begin{equation}\label{objective}
\Theta(w)\triangleq\sum_{i\in I} H_i(w_i)+\sum_{k\in K}\mbox{$\frac{I^k}{\alpha^k}$}\,\EU^k(\bar u^k_0,\bar u^k_1-\alpha^k w_1,,\ldots,\bar u^k_n-\alpha^k w_n)
\end{equation}
where ${\displaystyle H_i(w_i)=\int_{w^0_i}^{w_i}\!\!\theta_i^{-1}(z) \,dz}$. \\
\end{theorem}

The proof is presented in Appendix A. We also present a refinement of Theorem \ref{thm1} for the case in which there is no possibility to opt out:\\

\begin{theorem}
Theorem \ref{thm1} remains valid when there is no opt-out alternative, 
under the following non-saturation condition:
\begin{equation}\label{nonsaturation}
\sum_{i\in I}\bar x_i>\sum_{k\in K}I^k.\\
\end{equation}
\end{theorem}

The proof is presented in Appendix A.\\

Hence, we have established that the general model, which assumes that choice probabilities follow random utility theory and that waiting times at facilities strictly increase with the corresponding total inflows, can be formulated as a strictly convex minimization problem. The model therefore has a unique solution describing the equilibrium choice probabilities and waiting times. The absence of alternative local minima reduces the interest to explore alternative (non-cooperative) game theoretical equilibria.\\

Let us now consider several specific elaborations to which these general results apply.

\noindent
{\bf Example 1 -- Multinomial Logit (MNL).} A popular discrete choice model, supported by Gnedenko's theorem on extremal value distributions, considers
i.i.d. Gumbel variables $\{\varepsilon^k_{i}:i\in I\}$  centered at 0 and with scale 
parameter $\beta_k$, in which case
\begin{eqnarray*}
\EU^{k}(u^k)&=&\mbox{$\frac{1}{\beta_k}\ln\left(\sum_{j=0}^n \exp(\beta_k u^k_{j})\right)$},\\
\frac{\partial\EU^{k}}{\partial u^k_{i}}&=&\frac{\exp(\beta_ku_{i}^k)}{\sum_{j=0}^n\exp(\beta_ku^k_{j})}. 
\label{mnlchoiceprobability}
\end{eqnarray*}
Notice that the latter equation gives the choice probabilities as per \eqref{gradprox}.\\
\vspace{2ex}

\noindent {\bf Example 2 -- Queuing models.} We can use any queuing model for the waiting time $w_i=\theta_i(x_i)$
as a function of the arrival rate $x_i$. The simplest one is an M/M/1 queue with $w_i=\frac{x_i}{\mu_i(\mu_i-x_i)}$.
Similarly, for an M/M/$s$ queue (as also analysed in \cite{zhang2019medical}) we have $w_i=\frac{1}{\mu_i}Q_s\!\left(\frac{x_i}{\mu_i}\right)$ where 
\vspace{-2ex}
$$Q_s(z)=\left(\sum_{k=0}^{s-1}\frac{z^k}{k!}+\frac{z^{s}}{s!}\frac{s}{(s-z)}\right)^{-1}\frac{z^{s}}{s!}\frac{s}{(s-z)^2}.$$
Our general result can also accommodate queues with finite capacity such as M/M/$s$/$K$ queues and 
queues with different arrival distributions, such as M/G/$s$ queues (see also \cite{marianov2008facility, kucukyazici2020incorporating} and the references therein). The only restriction is that $\theta_i(\cdot)$ must be strictly increasing and continuous.

\section{Case Study in Urban China}

\subsection{Model Specification}

The case considers hospital choice in urban China. In line with the health system challenges and preceding research, the case study focuses on choice probabilities and waiting times per facility level rather than per facility. The indices $i\in\{1,2,3\}$ of facility types therefore will represent the primary, secondary, and tertiary hospitals, and $i=0$ refers to opting out. 

The patient choice related health system challenges center around the perceived weakness of primary care facilities. The low utility attached to the primary care level institutions by China's urban population -- and to a lesser extent to the level of secondary hospitals \cite{liu2020impact} -- causes 1) under-utilization of primary care, 2) over-utilization of tertiary care, and 3) undesirable avoidance of care, i.e., high opt-out rates \cite{meng2019can, ta2020trends, li2017primary, li2020quality, liu2018factors}. This problematic situation has resisted large scale healthcare reforms which mainly aimed to change choice behavior by improving the skills level of the medical staff and the quality of the medical equipment at primary care facilities. 

According to the DCE study conducted by Liu et al., skills level of medical staff and quality of medical equipment are indeed the two main determinants of patient choice for the population of 24 million of Shanghai \cite{liu2020impact, chinagovstats}). 

This case study uses the MNL model, data, and results obtained by the aforementioned DCE study \cite{liu2020impact}. To capture the endogeneity between choice probability and waiting time it additionally includes an M/M/1 queuing model. 

The case study thus advances on the previous DCE study and policy analysis therein which disregarded the endogeneity and included fixed waiting times as choice determinants \cite{liu2020impact}. We thus aim to present a more accurate assessment of the policy interventions implemented to remedy the health system challenges and to present insight in the (lack of) effectiveness of these interventions. Moreover the case explores the possible effectiveness (efficacy) of alternative interventions. \\

\subsection{MNL Model, Parameters and Data }

The patient choice data were collected in a DCE study for which patients gave their consent and for which ethics approval was obtained from the Shanghai General Hospital Medical Ethical Review Committee (no 2017KY207) \cite{liu2020impact}. 

In the MNL model, the utility of a facility level is determined by  the waiting time and by the factors presented in the top half of Table \ref{tablehospitalfactordata} (\textit{cf.} \cite{liu2020impact}). 

\begin{table}[h]
\begin{center}
{\small
\begin{tabular}{|l|c|c|c|}
\hline
 & \textbf{Primary} & \textbf{Secondary} & \textbf{Tertiary} \\
 \hline\hline
\textbf{Facility Attributes} & & & \\
\textit{Facility Size} & Small & Medium & Large \\
\textit{Out of Pocket Cost (RMB)} & 59 & 88 & 105 \\
\textit{Typical Skills Level} & Junior & Senior & Expert \\
\textit{Status of Equipment} & Outdated & Standard & Advanced \\
\textit{Travel Time  (min.)} & 15 & 15 & 15 \\
\textit{Total Visit Time (hours)} & 1 & 3 & 5 \\ 
\hline\hline
\textbf{Parameters} &&&\\
\textit{Service Rate} (patients/hour) & 10 & 10 & 12\\
\textit{Visit Time Other Than Waiting (min.)} & 34 & 88 & 87 \\
\textit{Zero Wait Utility for Mild Patients} & 0.207 & 0.417 & -0.259 \\
\textit{Zero Wait Utility for Severe Patients} & -0.257 & 0.089 & 0.773 \\
\textit{Waiting Time Multipliers}& 3 & 5 & 7 \\
\textit{Number of Facilities} & 1009 & 105 & 47\\
\textit{Number of Doctors per Facility} & 6.76 & 40.8 & 98.9 \\ 
\textit{Doctor Time Allocated to First Visits (\%)} & 50 & 50 & 50 \\
\hline
\end{tabular}
}
\end{center}
\caption{Attributes and parameters of facilities for all three levels}
\label{tablehospitalfactordata}
\end{table}

Focusing on the interplay between choice probability and waiting time per intervention, we recast the utility functions to the form $u^k_i=\bar u^k_i-\alpha^k w_i$, in which the only explicit factor is the waiting time $w_i$ and the reference utility $ \bar u^k_i$ incorporates the utilities from all other factors. Thus, $u^k_i$ represents the utility patients of type $k \in K$ attach to facilities of level $i \in I$. 

The original DCE study in \cite{liu2020impact} considered the factor sojourn time which is composed of the waiting time plus the visit time. The latter is assumed not to be affected by waiting time and  therefore included in the reference utility $\bar u^k_i$. 

As the urban citizens of Shanghai typically live close to facilities of each level and travel time is not among the most important factors determining facility choice, we assume a generic travel time of 15 minutes and distinguish patient types only by the evidence based factor perceived severity \cite{liu2018patients, liu2020impact}. More specifically, we consider two patient types: patients perceiving mild disease and patients perceiving severe disease. The patient type dependent utilities that differ significantly for different parameter levels are presented in Table \ref{tablepatientdata}.

\begin{table}
\begin{center}
\begin{tabular}{|l|c|c|}
\hline
\textbf{Factor Utilities} & \textbf{Mild} & \textbf{Severe} \\ 
\hline\hline
\textit{Opt Out Utility} & 2.499 & -6.024 \\ 
\textit{Waiting Time Sensitivities} $\alpha^k$ & 0.232 & 0.0995  \\
\hline\hline
\textit{Junior Doctor Utility} & -0.277 & -0.05 \\
\textit{Senior Doctor Utility} & 0.199 & -0.089 \\
\textit{Expert Doctor Utility} & 0.078 & 0.139 \\ 
\hline\hline
\textit{Obsolete Equipment Utility} & -0.275 & -0.43 \\
\textit{Advanced Equipment Utility} & 0.275 & 0.43 \\
\hline
\end{tabular}
\end{center}
\caption{Data determining the utility per patient type}
\label{tablepatientdata}
\end{table}

The sensitivity to waiting time varies depending on the disease severity, \textit{i.e.}, across patient types $k \in K$, but is independent of facility level. Opt out utilities also differ per patient type. The patient type dependent waiting time sensitivities and opt out utilities are also included in \ref{tablepatientdata}.\\

Recalling \eqref{mnlchoiceprobability} from Example 1 above, the MNL choice probabilities $\pi_{i}^k$ can now be expressed as:
\[ \pi^k_i =  \frac{\exp(\beta_k(\bar u^k_i - \alpha^k w_i))} {\sum_{j=0}^{3} \exp(\beta_k(\bar u^k_j - \alpha^k w_j))}. \] 

\subsection{M/M/1 Model, Parameters and Data }

Focusing the case study on patients who seek care for a first visit and hence without appointment with a pre-assigned physician, we assume the assignment of patients to physicians happens as soon as the patient is registered, \textit{i.e.}, upon arrival at the facility. Moreover, we assume that this causes patients to be evenly distributed among the doctors for each level. Thus, each registered patient joins a single server queue. More specifically, we assume arrivals and service rates are exponentially distributed and thus adopt an M/M/1 model and consider alternative queuing models in the discussion. The data required for the M/M/1 model are presented in the bottom half of Table \ref{tablehospitalfactordata}, as further outlined below.

As stated in Example 2 above, the M/M/1 model yields that expected waiting time equals $w_i=\frac{x_i}{\mu_i(\mu_i-x_i)}$ where $x_i$ is the arrival rate of patients for a doctor at facility level $i \in I$ and $\mu_i$ their service rate. 

Data on the service rates $\mu_i$ were provided by four medical doctors, two from primary care facilities, one from a secondary hospital, and one from a tertiary hospital in Shanghai in June 2019. The obtained data indicate that maximum service rates are normally around 12 patients per hour for each of the three facility levels. As tertiary hospitals are especially known for long waits and short service times \cite{chao2017healthcare, liu2018factors}, the base service rates are set at 12 patients per hour for tertiary care facilities and 10 patients per hour for primary and secondary level facilities.  

The same medical doctors also provided data on the average number of on duty doctors per day and documentation on the total number of licensed medical doctors per facility. The numbers of facilities per level were obtained from \cite{shanghaistats2018}.

The service capacity per level can then be determined by the number of facilities, the number of physicians per facility, and the fraction of physician time spent on first visit consults, which are all presented in Table \ref{tablehospitalfactordata}. Moreover, we assume 261 working days of 8 hours per year.\\

The incidence of patients who perceive to have an illness and consider visiting a healthcare facility is taken from the national census \cite{chinagovstats2013}. The division of these patients over the two severity types is determined by solving \eqref{gradprox}-\eqref{wait} for the base case model and data with the additional requirement that the solution yields the evidence based reference waiting times of the DCE on which it was based \cite{liu2020impact}. This yielded that 47.9 percent of the 160,432,700
patients in Shanghai per year who consider to visit a healthcare facility or to opt out instead, perceive their disease as minor.

The parameter values for the queuing model presented in the bottom half of Table \ref{tablehospitalfactordata} additionally reflect that primary care consults typically require a follow up visit to the hospital pharmacy and the cashier, secondary hospitals typically perform diagnostic services such as imaging and lab tests, and tertiary hospitals have the most elaborate and time consuming pathways \cite{liu2018patients, li2017primary, zou2017analysis}. Correspondingly, we assume that patients will join additional queues for these follow-up services (\textit{e.g.} at the pharmacy or for the lab test) as indicated by the waiting time multipliers depicted in Table \ref{tablehospitalfactordata}. 

\subsection{Comparative Intervention Analysis Methods}

The case study evaluates the effects on choice probabilities and waiting times of policy interventions which are modeled through modifications of the (default) parameter values depicted in Tables \ref{tablehospitalfactordata}, \ref{tablepatientdata}. Thus, each intervention is modeled by specifying a set of values for the utilities of all choice determinants with the exception of (the endogenous) waiting time and updating the reference utilities $\bar u^k_i$ accordingly. 

The comparative analysis presents both the results from solving  \eqref{gradprox}-\eqref{wait} with choice probabilities  determined by an MNL model and waiting times determined by an M/M/1 queuing model as explained above, and the results obtained by using an MNL model only with the constant reference total visit times of 1,3, and 5 hours for primary, secondary and tertiary care respectively from the DCE study \cite{liu2020impact}.\\

We compare the results for the base scenario and five policy interventions. The first three of these interventions aim to strengthen primary care and model the substantial efforts made as part of the national health reform to upgrade the medical skills of the physicians and to upgrade the equipment of primary care facilities. Skill and equipment are the factors evidenced to impact choice probabilities most \cite{li2017primary, liu2020impact, healthreform2020}. 

\begin{itemize}
    \item {\bf Upskill:} Maximally upskill primary care physicians from the \textit{mostly junior} level to the \textit{mostly expert} level of tertiary care physicians. To avoid secondary care having less skilled physicians than primary care, the secondary care physicians are upskilled to the \textit{mostly expert} level as well. The required upskilling effort can be considered to be unrealistically ambitious. Hence, the relevance of the intervention results is not to model the effectiveness of actual interventions but rather to estimate the maximum effect attainable by these actively practiced upskilling interventions.
\item {\bf Upgrade:} Maximally upgrade primary care equipment from the level \textit{obsolete} to the \textit{advanced} level of tertiary care. To avoid secondary care having poorer equipment than primary care, the secondary care equipment is upgraded to \textit{advanced} as well. Again, the required upgrading effort can be viewed as unrealistic, \textit{e.g.}, for being too costly, and the relevance of the analysis is to provide insight into the maximum effectiveness of the on going equipment upgrading interventions.
 \item {\bf Upskill and Upgrade:} The third intervention combines the previous two and therefore gives insight into the maximum effect attainable by combined upskilling and upgrading.
\end{itemize}   

As both evidence from practice and the comparative analysis results below illustrate, the (maximum) effectiveness of upskilling and upgrading interventions is limited. A first main underlying cause appears to be the high utility patients perceiving a mild disease attach to opting out (self care). Patients perceiving a mild disease attach a utility difference between opting out and a zero wait primary care visit of 2.3, more than twice the utility gain from maximum upskilling and upgrading (0.903) and equivalent to 10 hours of waiting time reduction (0.232 per hour) (see Tables \ref{tablepatientdata} and \ref{tablehospitalfactordata}). As a result more than three out of four such patients opt out in expectation, and this choice probability is only weakly impacted by upskilling and upgrading. Intervention four therefore models a health promotion campaign that reduces the utility of self care (opting out) by half.  

A second important underlying cause is that utility is quite insensitive to waiting times for patients perceiving a severe disease. One hour of waiting time reduction reduces the utility of a facility by 0.1 (approximately), while the utility difference between primary and tertiary care is more than ten times larger for these patients. For realistic waiting times which correspond to being serviced on the day of visiting, the probability of choosing primary care therefore remains low and the probability to choose tertiary care remains high for these patients. 

One may argue that the utility of waiting should only depend on the person and not on the severity of the illness. Then, the current difference between the waiting time utilities for perceiving mild and severe disease are an artifact that may reflect other more lowly valued characteristics of primary and secondary care not already captured by the utility factors presented in Table \ref{tablehospitalfactordata}. Uniform waiting time sensitivities, which disregard perceived disease severity, then more realistically reflect waiting time sensitivities and may apply to the aspired future state of the health system in which it provides balanced access to care. 

Thus we have two more interventions:
\begin{itemize}
    \item {\bf Health Promotion Campaign:} This campaign discourages self care for patients perceiving a mild disease and halves the utility of opting out from 2.499 to 1.250. 
\item {\bf Uniform Waiting Time Sensitivity:} This intervention does not reflect a health reform but rather adjusts the model assuming that the current long waiting times may no longer be valid in the future and adjust the present low sensitivities to waiting time for patients perceiving severe disease to become equal to the waiting time sensitivity of patients perceiving mild disease, 0.232.
\end{itemize}   

Appendix B additionally presents analysis results for an intervention in which upskilling is only for primary care physicians and thus from \textit{mostly junior} to \textit{mostly senior}, which is the skills level for secondary care. Moreover, the same appendix presents results obtained when combining the upskilling and upgrading interventions with interventions four and five. \\

As the data on service use and hospital capacity are not based on strong evidence, we conducted a sensitivity analysis for each of the interventions studied. For each intervention study, we generated 1,000 problem instances in which, for each of the facility levels, the waiting time multipliers (demand) and fraction of on duty doctors per day (supply) experience a uniformly distributed random perturbation of between -10\% and 10\%.

This sensitivity analysis provides insights into the robustness of the results. Moreover, it allows to test whether the differences between the results obtained by solving \eqref{gradprox}-\eqref{wait} are significantly different from the solutions obtained when solving the MNL model only. We present results on two statistical tests for the significance of difference. First, we conducted a sign test, testing the null hypothesis that there are no differences in choice probabilities and waiting times between the two models. For each of these parameters and for 1000 instances, this hypothesis is rejected if the number of instances for which the difference is positive (negative) is at least 526 (at most 473) because this number has a probability of less than 0.05. These results are indicated with `+' and `-', in case of rejection. A stronger variant of the same test rejects the hypothesis if the number of instances for which the difference is positive (negative) is at least 537 (at most 462), which has a probability of less than 0.01 under the null hypothesis. 

Even though the service use and capacity perturbations are sampled from uniform distributions which are symmetric around the mean values of the base scenario, the effects on waiting time can be highly nonlinear and asymmetric, if only because M/M/1 waiting time is superlinear in the arrival rate per doctor ($w_i=\frac{x_i}{\mu_i(\mu_i-x_i)}$). Thus, a second test rejects the null hypothesis of zero difference for choice probability or waiting time, if  the value is outside of a 95 \% confidence interval around the mean of the variable under consideration. This occurs when at most 25 of the 1000 problem instances report a difference below zero or above zero. Notice that this second test is much more demanding. Both tests are robust as they make no assumptions about the distributions of the (differences between) outcome parameters.

The Python code used for the computational studies is available on Github \cite{pythoncode}.

\subsection{Comparative Intervention Analysis Results}

Figure \ref{baselineresults} presents the results of the base case with the originally collected data and the results of 1000 perturbed model instances. Columns two and three present the choice probabilities and waiting times for the MNL model and for model \eqref{gradprox}-\eqref{wait} respectively for the unperturbed base instance. These columns present the choice probabilities $\pi_i^k=\mathbb{P}(i|k)$ denoted in column 1, starting with the probability of opting out by patients perceiving mild disease $\mathbb{P}(OO|M)$, to primary care for patients perceiving mild disease $\mathbb{P}(1|M)$, et cetera, to tertiary care for patients perceiving severe disease $\mathbb{P}(3|S)$. These probabilities confirm that the calibrated base instance of model \eqref{gradprox}-\eqref{wait} accurately models the empirical data on which the DCE study is based. This is further confirmed by the virtually equal waiting times reported by both models for each of the three hospital levels in the corresponding bottom rows of columns two and three.

\begin{table}
\centering
{\scriptsize
\begin{tabular}{|l|c|c|c|c|c|c|}\hline
\textbf{Baseline}&MNL&1-4&$\mu$\,MNL&$\mu$\,1-4&Sign Test&Nonzero\\ \hline\hline
$\mathbb{P}(OO|M)$&0.7757&0.7759&0.7757&0.7740&{\tiny-\,-}&\\
$\mathbb{P}(1|M)$&0.0784&0.0784&0.0784&0.0783&{\tiny-\,-}&\\
$\mathbb{P}(2|M)$&0.0967&0.0966&0.0967&0.0972&{\tiny +\,+}&\\
$\mathbb{P}(3|M)$&0.0492&0.0491&0.0492&0.0505&{\tiny+\,+}&\\\hline
\hline
$\mathbb{P}(OO|S)$&0.0006&0.0006&0.0006&0.0006&{\tiny-\,-}&\\
$\mathbb{P}(1|S)$&0.1917&0.1918&0.1917&0.1906&{\tiny-\,-}&\\
$\mathbb{P}(2|S)$&0.2709&0.2709&0.2709&0.2701&{\tiny -\,-}&\\
$\mathbb{P}(3|S)$&0.5368&0.5367&0.5368&0.5387&{\tiny+\,+}&\\\hline
\hline
$W(1)$&0.43&0.43&0.43&0.43&{\tiny-\,-}&\\
$W(2)$&1.53&1.52&1.50&1.48&{\tiny-\,-}&\\
$W(3)$&3.54&3.54&3.58&3.43&{\tiny +\,+}&\\\hline
\end{tabular}
}
\caption{Baseline choice probabilities and waiting times}
\label{baselineresults}
\end{table}

Columns four and five of Table \ref{baselineresults} present the average values over 1000 randomly perturbed problem instances for the MNL model and model \eqref{gradprox}-\eqref{wait} respectively. The choice probabilities in column four are identical to the choice probabilities in the second column as the MNL model does not consider effects of any changes in waiting time on the choice probabilities. The bottom rows show that the average waiting times differ slightly from the unperturbed base instance. However, one might argue that the column two and four waiting times should be disregarded as the MNL model assumes there is no change in the utility from waiting times and hence no change in waiting times. Thus, below we focus on the choice probabilities in the model comparison. Column five shows that the choice probabilities and waiting times differ slightly from those of column three when using model \eqref{gradprox}-\eqref{wait}.

Column six of Table \ref{baselineresults} shows the result of the sign test, where  `++' and`$--$' mean that the null hypothesis of zero difference holds with probability less than 0.01 and single `+', `-' indicate the same yet with probability at most 0.05. Thus `(+)+'s indicate that \eqref{gradprox}-\eqref{wait} results are significantly higher than MNL results and `(-)-'s indicate that these are significantly smaller. Empty cells in column six indicate that there are no significant differences. Even for this base case for which the model \eqref{gradprox}-\eqref{wait} is calibrated to reproduce the MNL outcomes and evidence-based DCE assumptions, the perturbations lead to differences in outcomes between the two models that cause the sign test to reject that hypothesis of equal outcomes ($p \leq 0.01$) for all choice probabilities and waiting times.

An `X' in column seven of Table \ref{baselineresults} indicates that the second test described in the methods section is rejected and any empty cell indicates it is not rejected, which is the case for all outcome variables in this base scenario.

Next, we turn to considering the model results for the upskilling and upgrading interventions. These results are presented in Table \ref{upskillingandupgradingresults}.

\begin{table}
    \centering
{\scriptsize
\begin{tabular}{|l|c|c|c|c|c|c|}\hline
\textbf{Upskill}&MNL&1-4&$\mu$\,MNL&$\mu$\,1-4&Sign Test&Nonzero\\ \hline\hline
$\mathbb{P}(OO|M)$&0.7587&0.7545&0.7587&0.7526&{\tiny-\,-}&\\
$\mathbb{P}(1|M)$&0.1094&0.1080&0.1094&0.1079&{\tiny-\,-}&\\
$\mathbb{P}(2|M)$&0.0838&0.0814&0.0838&0.0822&{\tiny -\,-}&\\
$\mathbb{P}(3|M)$&0.0481&0.0560&0.0481&0.0573&{\tiny+\,+}&$\times$\\\hline
\hline
$\mathbb{P}(OO|S)$&0.0005&0.0005&0.0005&0.0005&{\tiny-\,-}&\\
$\mathbb{P}(1|S)$&0.2087&0.2020&0.2087&0.2009&{\tiny-\,-}&$\times$\\
$\mathbb{P}(2|S)$&0.3068&0.2949&0.3068&0.2944&{\tiny -\,-}&$\times$\\
$\mathbb{P}(3|S)$&0.4840&0.5026&0.4840&0.5042&{\tiny+\,+}&$\times$\\\hline
\hline
$W(1)$&0.47&0.46&0.47&0.46&{\tiny-\,-}&$\times$\\
$W(2)$&1.77&1.62&1.72&1.57&{\tiny-\,-}&$\times$\\
$W(3)$&2.41&2.85&2.36&2.76&{\tiny +\,+}&$\times$\\\hline
\multicolumn{7}{}{}\\[2ex]\hline
\textbf{Upgrade}&MNL&1-4&$\mu$\,MNL&$\mu$\,1-4&Sign Test&Nonzero\\ \hline\hline
$\mathbb{P}(OO|M)$&0.7129&0.7046&0.7129&0.7029&{\tiny-\,-}&$\times$\\
$\mathbb{P}(1|M)$&0.1249&0.1194&0.1249&0.1193&{\tiny-\,-}&$\times$\\
$\mathbb{P}(2|M)$&0.1170&0.1110&0.1170&0.1110&{\tiny -\,-}&\\
$\mathbb{P}(3|M)$&0.0452&0.0650&0.0452&0.0657&{\tiny+\,+}&$\times$\\\hline
\hline
$\mathbb{P}(OO|S)$&0.0004&0.0004&0.0004&0.0004&{\tiny-\,-}&$\times$\\
$\mathbb{P}(1|S)$&0.3220&0.3003&0.3220&0.2994&{\tiny-\,-}&$\times$\\
$\mathbb{P}(2|S)$&0.2960&0.2752&0.2960&0.2755&{\tiny -\,-}&$\times$\\
$\mathbb{P}(3|S)$&0.3816&0.4241&0.3816&0.4247&{\tiny+\,+}&$\times$\\\hline
\hline
$W(1)$&0.62&0.58&0.61&0.57&{\tiny-\,-}&$\times$\\
$W(2)$&2.03&1.69&1.98&1.65&{\tiny-}&$\times$\\
$W(3)$&1.49&1.91&1.46&1.86&{\tiny +\,+}&$\times$\\\hline
\multicolumn{7}{}{}\\[2ex]\hline
\textbf{Upskill \&} &&&&&&\\
\textbf{Upgrade} &MNL&1-4&$\mu$\,MNL&$\mu$\,1-4&Sign Test&Nonzero\\\hline\hline
$\mathbb{P}(OO|M)$&0.6856&0.6791 &0.6856 & 0.6776&{\tiny-\,-}&$\times$\\
$\mathbb{P}(1|M)$&0.1713&0.1611 &0.1713 &0.1611&{\tiny-\,-}&$\times$\\
$\mathbb{P}(2|M)$&0.0997&0.0927 &0.0997 &0.0938&{\tiny -\,-}&$\times$\\
$\mathbb{P}(3|M)$&0.0435&0.067 &0.0435 1&0.0676&{\tiny+\,+}&$\times$\\\hline
\hline
$\mathbb{P}(OO|S)$&0.0003&0.0004&0.0004&0.0004&{\tiny-\,-}&$\times$\\
$\mathbb{P}(1|S)$&0.3403 &0.3158 &0.3403 &0.3151&{\tiny-\,-}&$\times$\\
$\mathbb{P}(2|S)$&0.3254&0.3005 &0.3254 &0.3011&{\tiny -\,-}&$\times$\\
$\mathbb{P}(3|S)$&0.3339 &0.3833 &0.3339 & 0.3834&{\tiny+\,+}&$\times$\\\hline
\hline
$W(1)$&0.72&0.66 &0.71 &0.65&{\tiny-\,-}&$\times$\\
$W(2)$&2.26&1.79 &2.20 &1.73&{\tiny-}&$\times$\\
$W(3)$&1.26 &1.61 &1.24 &1.50&{\tiny +\,+}&$\times$\\\hline
\end{tabular}

}

\caption{Choice probabilities and waiting times after upskilling and upgrading interventions}
\label{upskillingandupgradingresults}
\end{table}


The most notable policy observation from these three modeled intervention studies might be that even full upskilling and upgrading only reduces the opt  out probability of mild patients from 0.78 to 0.67 (or to 0.68 according to the MNL model which ignores waiting time effects). However, the relative changes in choice probabilities among the hospital levels tend to be larger. For instance, the original probability of patients perceiving severe disease to choose tertiary care was 0.54 and decreases to 0.50 (upskill), to 0.42 (upgrade), and to 0.38 (upskill and upgrade) for the successive interventions studied. The MNL model finds this probability to be as low as 0.33 and both significance tests reveal that this reduction appears significantly overestimated.

More generally we can observe that differences between the policy effects obtained by the models are always significant according to the sign test, and significant for the far majority on the second test as well. As a rule, the MNL model significantly overestimates the policy effects. There are some irregularities, for instance, because the utility severe patients attach to physicians is not increasing with the skills level. They value senior doctors the least. The results also indicate that the combination of full upskilling and upgrading goes a long way to balance the flow of patients perceiving severe disease but is hardly effective to prevent patients perceiving mild disease from opting out.\\

The Health Promotion intervention models the reduction of the opt out utility for patients perceiving mild disease. The modeled intervention results are presented in Table \ref{twomoreinterventionresults}. Columns three and five show the effectiveness of this intervention to reduce opting out among patients perceiving mild disease from 0.78 to around 0.54. The increased opting in leads to considerable increases in waiting times, especially for secondary and tertiary for which patients arrival rates were already close to facility level capacity. 
 
 For this health promotion intervention, the limitations of the MNL model again cause an overestimation of the effectiveness. Moreover,
Table \ref{twomoreinterventionresults} reveals another main shortcoming of the MNL model. It may yield choice probabilities that cause patient flows above capacity and thus unstable queues with infinite waiting times. Column four shows that only 739 out of the 1000 perturbed instances had feasible solutions. In fact, the feasibility limit applied is slightly more restrictive than requiring finite waiting times. All instances with a waiting time of more than ten hours are considered as infeasible. Such because regular hospital opening hours are limited to 10 hours and longer waiting might imply waiting until the next day. The estimated sensitivities of the DCE are unlikely to be valid for such situations, despite the fact that such exceptionally long waiting has occurred in Chinese urban health systems in the recent past. Appendix B shows that when combining the interventions of upgrading equipment and health promotion, even the undisturbed results (column two) are infeasible.

Finally, the same table shows the effects of changing the waiting time sensitivities of patients perceiving severe disease to equalize the sensitivities of mild patients. These effects are comparable to the effect of complete upskilling. Hence, they are relatively modest, with the exception of a substantial reduction in the choice probabilities of tertiary care by patients perceiving severe disease and the corresponding tertiary care waiting times. Because the MNL model assumes waiting times to be fixed, it fails to register any effects of modified waiting time sensitivities on choice probabilities. 

For both these additional interventions, the sign test qualifies all differences between the outcomes of the two models as highly  significant, as is mostly confirmed by the stricter second test.

\begin{table}
    \centering
{\scriptsize
\begin{tabular}{|l|c|c|c|c|c|c|}\hline
\textbf{Health}&&&&&&\\
\textbf{Promotion}&MNL&1-4&$\mu$\,MNL&$\mu$\,1-4&Sign Test&Nonzero\\ \hline\hline
$\mathbb{P}(OO|M)$&0.4980&0.5405&0.4980&0.5348&{\tiny+\,+}&$\times$\\
$\mathbb{P}(1|M)$&0.1755&0.1859&0.1755&0.1844&{\tiny+\,+}&$\times$\\
$\mathbb{P}(2|M)$&0.2164&0.1818&0.2164&0.1851&{\tiny -\,-}&$\times$\\
$\mathbb{P}(3|M)$&0.1101&0.0917&0.1101&0.0957&{\tiny-\,-}&\\\hline
\hline
$\mathbb{P}(OO|S)$&0.0006&0.0006&0.0006&0.0006&{\tiny+\,+}&$\times$\\
$\mathbb{P}(1|S)$&0.1917&0.2079&0.1917&0.2053&{\tiny+\,+}&$\times$\\
$\mathbb{P}(2|S)$&0.2709&0.2661&0.2709&0.2656&{\tiny -\,-}&\\
$\mathbb{P}(3|S)$&0.5368&0.5254&0.5368&0.5285&{\tiny-\,-}&\\\hline
\hline
$W(1)$&0.51&0.54&0.51(739)&0.53&{\tiny+\,+}&$\times$\\
$W(2)$&4.08&2.62&4.19(739)&2.50&{\tiny-\,-}&$\times$\\
$W(3)$&6.92&4.66&5.24(739)&4.46&{\tiny -\,-}&\\\hline
\multicolumn{7}{}{}\\[2ex]\hline
\textbf{Uniform}&&&&&&\\
\textbf{Waiting Time}&&&&&&\\
\textbf{Sensitivities}&MNL&1-4&$\mu$\,MNL&$\mu$\,1-4&Sign Test&Nonzero\\ \hline\hline
$\mathbb{P}(OO|M)$&0.7757&0.7697&0.7757&0.7678&{\tiny-\,-}&$\times$\\
$\mathbb{P}(1|M)$&0.0784&0.0773&0.0784&0.0773&{\tiny-\,-}&$\times$\\
$\mathbb{P}(2|M)$&0.0967&0.0951&0.0967&0.0961&{\tiny -\,-}&\\
$\mathbb{P}(3|M)$&0.0492&0.0577&0.0492&0.0588&{\tiny+\,+}&$\times$\\\hline
\hline
$\mathbb{P}(OO|S)$&0.0006&0.0007&0.0006&0.0007&{\tiny+\,+}&$\times$\\
$\mathbb{P}(1|S)$&0.1917&0.2252&0.1917&0.2227&{\tiny+\,+}&$\times$\\
$\mathbb{P}(2|S)$&0.2709&0.2753&0.2709&0.2750&{\tiny +\,+}&\\
$\mathbb{P}(3|S)$&0.5368&0.4986&0.5368&0.5015&{\tiny-\,-}&$\times$\\\hline
\hline
$W(1)$&0.43&0.46&0.43&0.45&{\tiny+\,+}&$\times$\\
$W(2)$&1.52&1.55&1.48&1.50&{\tiny+\,+}&\\
$W(3)$&3.54&2.80&3.55&2.72&{\tiny -\,-}&\\\hline
\end{tabular}
}

    \caption{Choice probabilities and waiting times after intervening on opt out utility and waiting time sensitivity}
    \label{twomoreinterventionresults}
\end{table}


\section{Discussion and Conclusion}

We discuss the results in the reverse order of presenting them, first interpreting case study results from urban China and then zooming out to the presented models, as well as present and future theoretical contributions on the relationship between patient choice and waiting time.

The inclusion of the mutual dependence between choice probabilities in the modeling of health policy intervention effects strongly suggests that the ongoing efforts to upskill doctors and upgrade equipment can contribute to alleviating the persistent problems of underutilization of primary care and overutilization of tertiary care. At the same time, the analysis reveals that these interventions likely fall short of resolving these problems.  In comparison to the previous MNL based policy analyses, the increase in utility of primary and secondary care facilities -- and thus the probability to choose facilities at these levels instead of opting out -- that results from these interventions are significantly smaller when accounting for endogeneity between waiting time and choice probability by solving \eqref{gradprox}-\eqref{wait}. An intuitive interpretation is that the utility increases caused by the interventions yield higher choice probabilities, which in turn imply longer waiting times that partially undo the initial utility increases. 

While the opt in probabilities remain low for patients perceiving mild disease for both models, the differences in effects of the intervention to upgrade equipment and upskill medical doctors are highly significant and larger in relative terms. For example, the MNL model estimates their probability to attend tertiary care is 0.0435 while solving  \eqref{gradprox}-\eqref{wait} yields a 0.0241 higher probability. Thus, the relative difference is more than 50 percent. 

The patients perceiving severe disease almost all opt in and the equally significant differences between model outcomes of upgrading and upskilling now also translate into larger absolute differences in the choice probabilities. For example, the MNL model estimates their probability to attend tertiary care after upgrading and upskilling is 0.3339 while solving  \eqref{gradprox}-\eqref{wait} gives a probability of 0.3833. This difference translates to around 4 million first visits to tertiary care per year.

Based on the above observations, we analysed two alternative interventions. The first suggestion is a health promotion campaign that diminishes the utility attached to self care and opting out. Our analysis reveals that such a campaign has much more efficacy potential than upskilling and upgrading. Moreover, our findings suggest that MNL based analysis significantly errs on the resulting choice probability changes and that its disregard of endogeneity can even yield choice probabilities that are infeasible because they violate hospital capacity limits. 

The MNL model is incapable of evaluating interventions which cause changes in waiting time sensitivities as it assumes waiting times to be constant.  Hence, the newly developed model is a welcome advancement to provide insight into the effect of remedying the waiting time insensitivity reported for patients perceiving severe disease \cite{liu2020impact}. The analysis shows that the effects of interventions to adjust waiting time sensitivities of patients perceiving severe disease to those of patients perceiving mild disease, resulting in uniform patient dependent sensitivities, are comparable to those obtained by full upskilling of physicians. This intervention follows the universal health coverage principles of the WHO which reaffirm the importance of access to strong primary care facilities for all and referral to higher levels as needed  \cite{ghebreyesus2018primary}. Further research into the determinants of waiting time (in)sensitivity and interventions to resolve the insensitivity is called for. \\

The case study analysis and in particular the analysis comparing the results obtained by the newly developed model with the results obtained by the MNL model only shows that the newly presented model can significantly more accurately capture the effects of interventions aimed at altering patient choice and corresponding health system utilization and performance. This confirms the importance of capturing the endogeneous relationships in the model for intervention analysis \cite{train2009discrete, antonakis2010making}. 

Regarding the data collection for the case study at hand, one may wonder whether revealed preference data would have yielded more accurate results. The choice set of the DCE limited waiting times to 1, 3, and 5 hours when eliciting stated preferences. Revealed preferences can potentially include a wider variety of combinations of choices and waiting times, and more realistic ones based on data per hospital rather than per hospital level and for different moments in time. Notice that this also requires obtaining corresponding revealed opt out choices, which is not trivial. 

Another difficulty with revealed preference data is that the attributes equipment, skills level, out of pocket costs, and hospital size, will be highly correlated as these are regulated to vary among hospital levels but not between hospitals of the same level. This renders the data to be unsuited for the analysis of the policy interventions targeting the individual attributes considered in our study. Combined with earlier evidence that stated preferences can predict revealed choices with a high degree of confidence for healthcare choices, we believe that stated preferences elicited through appropriately designed (MNL) models and corresponding choice sets may be most valuable to collect data for policy intervention studies in waiting time and choice interactions \cite{de2020can}.

Regarding the choice of queuing model for the case study at hand, one may argue that the M/M/1 queuing model which assumes patients distribute uniformly over the physicians per level has limited validity and moreover limits generalization of the results to other settings. Indeed, service rates or arrival rates may not be exponentially distributed and the service regime may be different in other settings. Example 2 of the modeling section already discusses the alternatives of M/M/$s$, /M/M/$s$/$K$ and M/G/$s$ queuing models. These examples illustrate that the proposed mathematical model is much more general than the specific model of the case study. It accommodates all queuing regimes in which waiting times are increasing in the choice probabilities, in combination with a variety of random utility based choice models. Theorem 3 extends the model to contexts in which opting out is not permitted. 

The general model presented still offers a variety of research directions worthwhile pursuing to strengthen and extend it. A main area of improvement lies in the assumption that hospital utility is a linear function of expected waiting time. We are not aware of evidence supporting this assumption, which is for instance implied by the MNL model. There are several indications that this relationship may not be linear. For travel time, for instance, there is evidence that any travel time less than 30 minutes is acceptable \cite{varkevisser2010assessing}. More generally, it may appear incorrect that a waiting time increase from 15 minutes to 45 minutes has the same effect on utility as a waiting time increase from 5 hours and 15 minutes to 5 hours and 45 minutes. Thus generalizations of the proofs of Theorems 2 and 3 to such alternative models, possibly with multiple equilibrium solutions, are called for. 

On the same thread, the utility of waiting time likely not only depends on expected waiting time but also on other moments such as the variance, or on the probability of long waiting. This is especially important as waiting times may vary considerably in practice for reasons of variability in demand and capacity as modelled in the perturbed instances of the computational study. Thus, we call for studies which include uncertainty of waiting time and nonlinear relations between expected waiting time and utility, for example using prospect theory based value functions which have been recently explored in a variety of studies on patient preferences \cite{rouyard2018risk, stolk2022patients}. 

Among the interesting extensions of the general model, choices and queuing beyond the first visit are of interest, as well as models in which service tasks and times vary depending on patient type (severity). 

Another valuable direction for extension is formed by the modeling of the learning dynamics of patients. We have assumed that equilibrium choice and waiting time will materialize but policy interventions  intend to disturb the equilibrium, for instance, to reduce waiting times. To analyse the effectiveness of interventions, it is therefore also of interest to know how patients learn about waiting times and how newly learned waiting times translate into updated utilities and choices. This matter is especially challenging as the expected waiting time is not easy to observe. The queue length at a given moment in time may be observable, but this does not easily translate into an accurate assessment of expected queue length or subsequent expected waiting time. Moreover, patients may typically have few waiting time observations from personal visits.  Added to the fact that learning processes may be slow and not fully rational, the envisioned equilibrium effects of policy interventions may not be realized until much later, if at all \cite{erev1998predicting}. Further research into the underlying strategic learning dynamics of hospital choice and how they can be promoted is thus an important area of future research \cite{young2004strategic, hart2005adaptive}. This research can borrow from related research in transportation science (see for instance \cite{cominetti2010payoff}
and references therein). 


\newpage

\noindent {\bf \Large Appendix A: Proofs}

\begin{proof}[\bf Proof of Theorem 1]
Let $U_i=u_i+\varepsilon_i$. For each realization of the random variables $\varepsilon_i$, the expression 
$$\varphi(u,\varepsilon)=\max\{u_{1}+\varepsilon_{1},\ldots,u_{n}+\varepsilon_{n}\}$$ is a maximum of linear functions
of $u$, hence it is a convex polyhedral function. Therefore its expectation, being an average of these, is also convex
$$\EU(u)=\int_{\RR^{n}}\varphi(u,\varepsilon)\, d\PP(\varepsilon).$$
Now, let $e_i$ be the $i$-th canonical vector in $\RR^n$. Then, when $t\to 0$ with $t>0$ one can readily check that 
the differential quotients 
of $\varphi(\cdot,\varepsilon)$ converge  monotonically towards an indicator function, namely
$$\frac{\varphi(u+t e_i,\varepsilon)-\varphi(u,\varepsilon)}{t}\to 
\left\{\begin{array}{cl} 1& \mbox{if }u_i+\varepsilon_i\geq u_j+\varepsilon_j,\forall j\\
0&\mbox{otherwise},
\end{array}\right.$$
and therefore Lebesgue's monotone convergence theorem yields
$$\lim_{t\to 0^+}\frac{\EU(u+te_i)-\EU(u)}{t}=\int_{\RR^{n}}\mathbbm{1}_{\{u_i+\varepsilon_i\geq u_j+\varepsilon_j,\forall j\}}\, d\PP(\varepsilon)=\PP(U_i\geq U_j,\forall j).$$
Similarly, for $t\to 0$ with $t<0$ we have
$$\lim_{t\to 0^-}\frac{\EU(u+te_i)-\EU(u)}{t}=\int_{\RR^{n}}\mathbbm{1}_{\{u_i+\varepsilon_i> u_j+\varepsilon_j,\forall j\}}\, d\PP(\varepsilon)=\PP(U_i> U_j,\forall j).$$
Since the random variables $\varepsilon_i$ are assumed continuous these two limits coincide and then 
$$\frac{\partial \EU}{\partial u_i}(u)=\PP(U_i\geq U_j,\forall j)=\PP(U_i> U_j,\forall j).$$
As every convex function which is Gateaux differentiable is automatically of class $C^1$ \cite{rockafellar1997convex}, this completes the proof.
\end{proof}

\begin{proof}[\bf Proof of Theorem 2]
Since the expected utilities $\EU^k(\cdot)$ are smooth and convex, while the functions $\theta_i^{-1}(\cdot)$ are strictly increasing,
it follows that $\Theta(\cdot)$ is strictly convex and smooth, with 
\begin{equation}\label{gradobj}
\frac{\partial\Theta}{\partial w_i}(w)=\theta_i^{-1}(w_i) -\sum_{k\in K}I^k\,\frac{\partial \EU^k}{\partial u^k_i}(\bar u^k_0,\bar u^k_1-\alpha^k w_1,\ldots,\bar u^k_n-\alpha^k w_n).
\end{equation}
Using \eqref{gradprox}-\eqref{utility} it follows that \eqref{reduced} is equivalent to $\nabla\Theta(w)=0$, and therefore 
the equilibria coincide with the minimizers of $\Theta(\cdot)$.

Since $\Theta(\cdot)$ is strictly convex it has at most one minimizer, and therefore
it remains to establish that the minimum is attained. This follows by showing that the 
recession function satisfies $\Theta^\infty(d)>0$ for all $d\in\RR^n\setminus\{0\}$.
Indeed, let us take $d\neq 0$. From Lebesgue's monotone convergence theorem, we have
$$(\EU^k)^\infty(d)=\max\{0,-\alpha^k d_1,\ldots,-\alpha^k d_n\}=-\alpha^k \min\{0,d_1,\ldots,d_n\},$$
and then using standard rules for computing recession functions we get
$$\Theta^\infty(d)=\sum_{i\in I}H_i^\infty(d_i) -\sum_{k\in K}I^k\min\{0,d_1,\ldots,d_k\}$$
with
$$H_i^\infty(d_i)=\lim_{t\to\infty}\theta_i^{-1}(td_i)d_i=
\left\{\begin{array}{cl}
\bar x_i d_i&\mbox{if }d_i>0\\
0&\mbox{if }d_i=0\\
+\infty&\mbox{if }d_i<0.
\end{array}\right.
$$
It follows that $\Theta^\infty(d)=+\infty$ whenever $d_i<0$ for some $i\in I$.
Otherwise, when $d_i\geq 0$ for all $i\in I$, we get $\Theta^\infty(d)=\sum_{i\in I}\bar x_i\,d_i>0$
which follows since $\bar x_i>0$ for all $i\in I$ and $d\neq 0$. This completes the proof. 
\end{proof}

\begin{proof}[\bf Proof of Theorem 3]
After removing the opt-out facility $i=0$ from the expression of the expected utilities
$$\EU^{k}(u^k)=\EE(\max\{u^k_{1}+\varepsilon^k_{1},\ldots,u^k_{n}+\varepsilon^k_{n}\}),$$
the characterization of equilibria as the unique minimizer of the strictly convex smooth 
function $w\mapsto\Theta(w)$ remains valid with the same proof. 

For the existence of equilibria, we observe that we still have $\Theta^\infty(d)=+\infty$ 
when $d_i<0$ for some $i\in I$, and therefore we only need to check that
\begin{equation}\label{nonsaturation2}
\sum_{i\in I}\bar x_i\,d_i>\sum_{k\in K}I^k\min\{d_1,\ldots,d_n\}\qquad(\forall \,d\geq 0, d\neq 0).
\end{equation}
Let $m=\min\{d_1,\ldots,d_n\}$.
The inequality \eqref{nonsaturation2} is again trivial when $m=0$, whereas for $m>0$
we have that the expression on the left is minimal when $d_i=m$ for all $i\in I$, 
in which case \eqref{nonsaturation2} reduces precisely to  \eqref{nonsaturation}.
\end{proof} 

\newpage

\vspace*{6cm}
{\bf \Large Appendix B: Intervention Study Results}

\includepdf[pages=-]{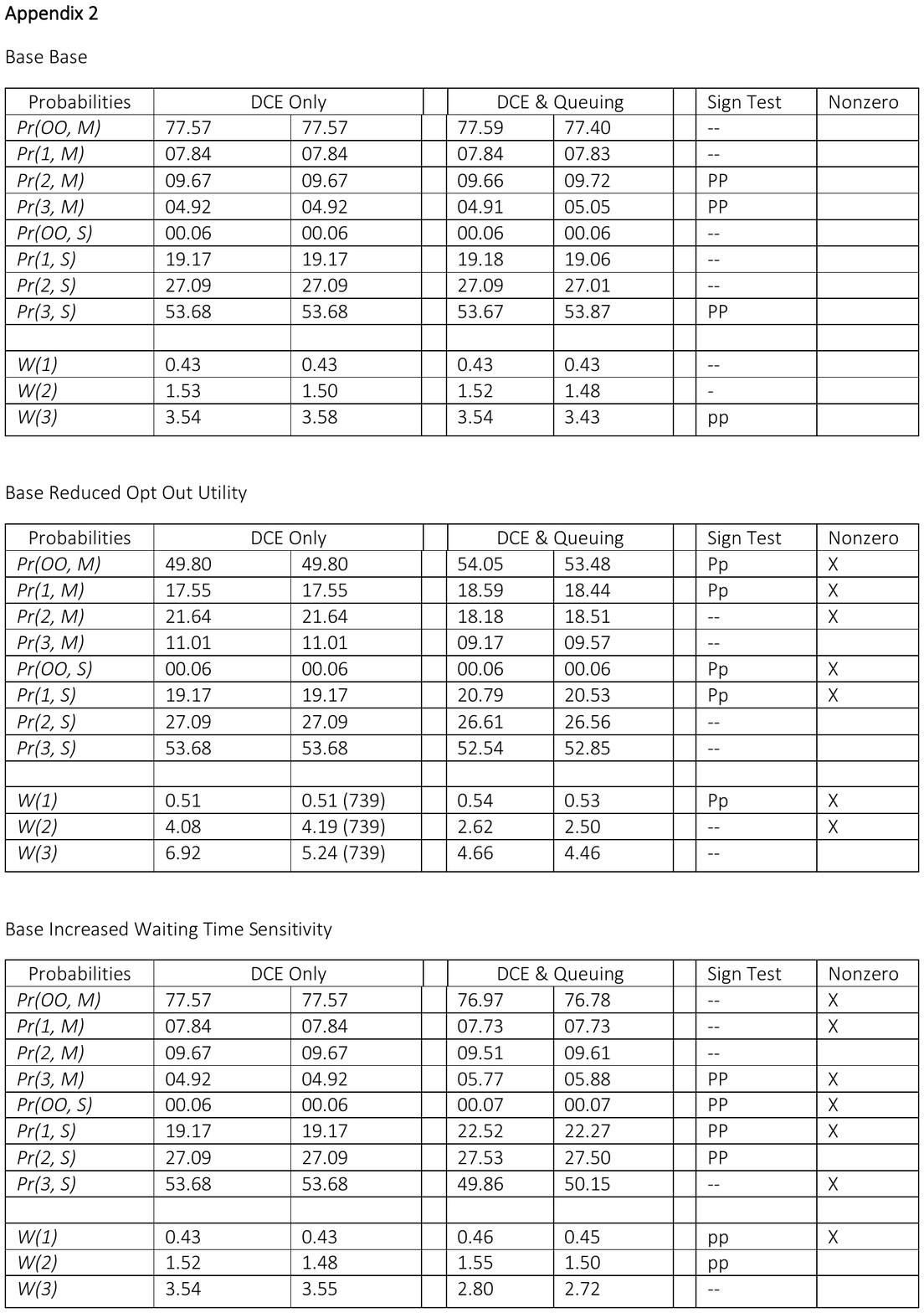}


{\small 
\begin{thebibliography}{10}

\bibitem{antonakis2010making}
John Antonakis, Samuel Bendahan, Philippe Jacquart, and Rafael Lalive.
\newblock On making causal claims: A review and recommendations.
\newblock {\em The leadership quarterly}, 21(6):1086--1120, 2010.

\bibitem{arsenault2020hospital}
Catherine Arsenault, Min~Kyung Kim, Amit Aryal, Adama Faye, Jean~Paul Joseph,
  Munir Kassa, Tizta~Tilahun Degfie, Talhiya Yahya, and Margaret~E Kruk.
\newblock Hospital-provision of essential primary care in 56 countries:
  determinants and quality.
\newblock {\em Bulletin of the World Health Organization}, 98(11):735, 2020.

\bibitem{berhane2015patients}
Adugnaw Berhane and Fikre Enquselassie.
\newblock Patients' preferences for attributes related to health care services
  at hospitals in amhara region, northern ethiopia: a discrete choice
  experiment.
\newblock {\em Patient preference and adherence}, 9:1293, 2015.

\bibitem{brown2015hospital}
Paul Brown, Laura Panattoni, Linda Cameron, Stephanie Knox, Toni Ashton, Tim
  Tenbensel, and John Windsor.
\newblock Hospital sector choice and support for public hospital care in new
  zealand: Results from a labeled discrete choice survey.
\newblock {\em Journal of health economics}, 43:118--127, 2015.

\bibitem{chao2017healthcare}
Jianqian Chao, Boyang Lu, Hua Zhang, Liguo Zhu, Hui Jin, and Pei Liu.
\newblock Healthcare system responsiveness in jiangsu province, china.
\newblock {\em BMC health services research}, 17:1--7, 2017.

\bibitem{cominetti2010payoff}
Roberto Cominetti, Emerson Melo, and Sylvain Sorin.
\newblock A payoff-based learning procedure and its application to traffic
  games.
\newblock {\em Games and Economic Behavior}, 70(1):71--83, 2010.

\bibitem{pythoncode}
Roberto Cominetti and Joris Van~de Klundert.
\newblock Python code for equilibrium computation.
\newblock
  \url{https://github.com/JorisAtWork/HospitalChoiceandWaitingTimeEquilbria/blob/main/Copy_of_HospitalChoice.ipynb}.

\bibitem{de2020can}
Esther~W de~Bekker-Grob, Bas Donkers, Michiel~CJ Bliemer, Jorien Veldwijk, and
  Joffre~D Swait.
\newblock Can healthcare choice be predicted using stated preference data?
\newblock {\em Social science \& medicine}, 246:112736, 2020.

\bibitem{domencich1975urban}
Thomas~A Domencich and Daniel McFadden.
\newblock {\em Urban travel demand-a behavioral analysis}.
\newblock North-Holland Publishing Company Limited, Oxford, England, 1975.

\bibitem{erev1998predicting}
Ido Erev and Alvin~E Roth.
\newblock Predicting how people play games: Reinforcement learning in
  experimental games with unique, mixed strategy equilibria.
\newblock {\em American economic review}, pages 848--881, 1998.

\bibitem{fischer2015understanding}
Sophia Fischer, Stefanie Pelka, and Ren{\'e} Riedl.
\newblock Understanding patients' decision-making strategies in hospital
  choice: Literature review and a call for experimental research.
\newblock {\em Cogent Psychology}, 2(1):1116758, 2015.

\bibitem{ghebreyesus2018primary}
Tedros~Adhanom Ghebreyesus, Henrietta Fore, Yelzhan Birtanov, and Zsuzsanna
  Jakab.
\newblock Primary health care for the 21st century, universal health coverage,
  and the sustainable development goals.
\newblock {\em The Lancet}, 392(10156):1371--1372, 2018.

\bibitem{hart2005adaptive}
Sergiu Hart.
\newblock Adaptive heuristics.
\newblock {\em Econometrica}, 73(5):1401--1430, 2005.

\bibitem{jones2022association}
Simon Jones, Chris Moulton, Simon Swift, Paul Molyneux, Steve Black, Neil
  Mason, Richard Oakley, and Clifford Mann.
\newblock Association between delays to patient admission from the emergency
  department and all-cause 30-day mortality.
\newblock {\em Emergency Medicine Journal}, 39(3):168--173, 2022.

\bibitem{kucukyazici2020incorporating}
Beste Kucukyazici, Yue Zhang, Amir Ardestani-Jaafari, and Lijie Song.
\newblock Incorporating patient preferences in the design and operation of
  cancer screening facility networks.
\newblock {\em European Journal of Operational Research}, 287(2):616--632,
  2020.

\bibitem{chinagovstats2013}
Sheng Laiyun et~al.
\newblock {\em An Analysis Report of National Health Services Survey in China}.
\newblock Center for Health Statistics and Information, 2013.

\bibitem{lee1985equilibrium}
Hau~Leung Lee and Morris~A Cohen.
\newblock Equilibrium analysis of disaggregate facility choice systems subject
  to congestion-elastic demand.
\newblock {\em Operations research}, 33(2):293--311, 1985.

\bibitem{li2020quality}
Xi~Li, Harlan~M Krumholz, Winnie Yip, Kar~Keung Cheng, Jan De~Maeseneer,
  Qingyue Meng, Elias Mossialos, Chuang Li, Jiapeng Lu, Meng Su, et~al.
\newblock Quality of primary health care in china: challenges and
  recommendations.
\newblock {\em The Lancet}, 395(10239):1802--1812, 2020.

\bibitem{li2017primary}
Xi~Li, Jiapeng Lu, Shuang Hu, KK~Cheng, Jan De~Maeseneer, Qingyue Meng, Elias
  Mossialos, Dong~Roman Xu, Winnie Yip, Hongzhao Zhang, et~al.
\newblock The primary health-care system in china.
\newblock {\em The Lancet}, 390(10112):2584--2594, 2017.

\bibitem{liu2020impact}
Yun Liu, Qingxia Kong, Shan Wang, Liwei Zhong, and Joris van~de Klundert.
\newblock The impact of hospital attributes on patient choice for first visit:
  evidence from a discrete choice experiment in shanghai, china.
\newblock {\em Health policy and planning}, 35(3):267--278, 2020.

\bibitem{liu2018factors}
Yun Liu, Qingxia Kong, Shasha Yuan, and Joris van~de Klundert.
\newblock Factors influencing choice of health system access level in china: a
  systematic review.
\newblock {\em PLoS One}, 13(8):e0201887, 2018.

\bibitem{liu2018patients}
Yun Liu, Liwei Zhong, Shasha Yuan, and Joris van~de Klundert.
\newblock Why patients prefer high-level healthcare facilities: a qualitative
  study using focus groups in rural and urban china.
\newblock {\em BMJ global health}, 3(5):e000854, 2018.

\bibitem{marianov2008facility}
Vladimir Marianov, Miguel R{\'\i}os, and Manuel~Jos{\'e} Icaza.
\newblock Facility location for market capture when users rank facilities by
  shorter travel and waiting times.
\newblock {\em European Journal of Operational Research}, 191(1):32--44, 2008.

\bibitem{martinez2019prolonged}
Diego~A Martinez, Haoxiang Zhang, Magdalena Bastias, Felipe Feijoo, Jeremiah
  Hinson, Rodrigo Martinez, Jocelyn Dunstan, Scott Levin, and Diana Prieto.
\newblock Prolonged wait time is associated with increased mortality for
  chilean waiting list patients with non-prioritized conditions.
\newblock {\em BMC public health}, 19(1):1--11, 2019.

\bibitem{mcintyre2020waiting}
Daniel McIntyre and Clara~K Chow.
\newblock Waiting time as an indicator for health services under strain: a
  narrative review.
\newblock {\em INQUIRY: The Journal of Health Care Organization, Provision, and
  Financing}, 57:0046958020910305, 2020.

\bibitem{meng2019can}
Qingyue Meng, Anne Mills, Longde Wang, and Qide Han.
\newblock What can we learn from china's health system reform?
\newblock {\em bmj}, 365, 2019.

\bibitem{chinagovstats}
Natural~Bureau of~Statistics~of China.
\newblock {\em China Statistical Yearbook 2019}.
\newblock China Statistics Press, 2019.

\bibitem{healthreform2020}
National Health~Commission of~the People's Republic~of China.
\newblock Notice of the national health commission on comprehensively promoting
  the construction of community health centers, 2020.

\bibitem{rockafellar1997convex}
R~Tyrrell Rockafellar.
\newblock {\em Convex analysis}, volume~11.
\newblock Princeton university press, 1997.

\bibitem{rouyard2018risk}
Thomas Rouyard, Arthur Attema, Richard Baskerville, Jos{\'e} Leal, and Alastair
  Gray.
\newblock Risk attitudes of people with manageable chronic disease: an analysis
  under prospect theory.
\newblock {\em Social Science \& Medicine}, 214:144--153, 2018.

\bibitem{shanghaistats2018}
MSB Shanghai.
\newblock {\em Shanghai statistical yearbook}.
\newblock China Statistics Press, 2018.

\bibitem{sheffi1978another}
Yosef Sheffi and Carlos~F Daganzo.
\newblock Another paradox of traffic flow.
\newblock {\em Transportation Research}, 12(1):43--46, 1978.

\bibitem{smith2018patient}
Honora Smith, Christine Currie, Pornpimol Chaiwuttisak, and Andreas Kyprianou.
\newblock Patient choice modelling: how do patients choose their hospitals?
\newblock {\em Health care management science}, 21(2):259--268, 2018.

\bibitem{stolk2022patients}
Aline~C Stolk-Vos, Arthur~E Attema, Michele Manzulli, and Joris~J van~de
  Klundert.
\newblock Do patients and other stakeholders value health service quality
  equally? a prospect theory based choice experiment in cataract care.
\newblock {\em Social Science \& Medicine}, 294:114730, 2022.

\bibitem{ta2020trends}
Yuqi Ta, Yishan Zhu, and Hongqiao Fu.
\newblock Trends in access to health services, financial protection and
  satisfaction between 2010 and 2016: Has china achieved the goals of its
  health system reform?
\newblock {\em Social Science \& Medicine}, 245:112715, 2020.

\bibitem{train2009discrete}
Kenneth~E Train.
\newblock {\em Discrete choice methods with simulation}.
\newblock Cambridge university press, 2009.

\bibitem{varkevisser2010assessing}
Marco Varkevisser, St{\'e}phanie~A van~der Geest, and Frederik~T Schut.
\newblock Assessing hospital competition when prices don't matter to patients:
  the use of time-elasticities.
\newblock {\em International Journal of Health Care Finance and Economics},
  10:43--60, 2010.

\bibitem{victoor2012determinants}
Aafke Victoor, Diana~MJ Delnoij, Roland~D Friele, and Jany~JDJM Rademakers.
\newblock Determinants of patient choice of healthcare providers: a scoping
  review.
\newblock {\em BMC health services research}, 12(1):1--16, 2012.

\bibitem{williams1977formation}
Huw~CWL Williams.
\newblock On the formation of travel demand models and economic evaluation
  measures of user benefit.
\newblock {\em Environment and planning A}, 9(3):285--344, 1977.

\bibitem{young2004strategic}
H~Peyton Young.
\newblock {\em Strategic learning and its limits}.
\newblock OUP Oxford, 2004.

\bibitem{zhang2019medical}
Yue Zhang and Derek Atkins.
\newblock Medical facility network design: User-choice and system-optimal
  models.
\newblock {\em European Journal of Operational Research}, 273(1):305--319,
  2019.

\bibitem{zhu2019exploring}
Jingrong Zhu, Jinlin Li, Zengbo Zhang, Hao Li, and Lingfei Cai.
\newblock Exploring determinants of health provider choice and heterogeneity in
  preference among outpatients in beijing: a labelled discrete choice
  experiment.
\newblock {\em BMJ open}, 9(4):e023363, 2019.

\bibitem{zou2017analysis}
DD~ZOU, ZF~ZHANG, et~al.
\newblock Analysis on outpatient stay in tertiary public hospitals in shanghai
  and relative patient satisfaction degree.
\newblock {\em Chinese Health Resources}, 20(6):464--468, 2017.

\end{thebibliography}
}
\end{document}